\documentclass[11pt]{article}

\usepackage{amsmath,amssymb,amscd}

\newtheorem{theorem}{Theorem}
\newtheorem{lemma}[theorem]{Lemma}
\newtheorem{prop}[theorem]{Proposition}

\begin{document}

\title{A universal enveloping for $L_\infty$-algebras.}
\author{Vladimir Baranovsky}
\date{June 10, 2007}
\maketitle

\begin{abstract}
For any $L_\infty$-algebra $L$ we construct an 
$A_\infty$-algebra structure on the symmetric coalgebra
$Sym^*_c(L)$ and prove that this structure 
satifies properties generalizing those of
the usual universal enveloping algebra. 
We also obtain an invariant contracting homotopy one the cobar 
construction of a symmetric coalgebra, by relating it to 
the combinatorics of permutahedra and  semistandard Young tableaux.  
\end{abstract}

\section{Introduction}

The purpose of this article is to generalize the universal enveloping
from Lie to $L_\infty$-algebras. One candidate is well-known: the cobar construction $\Omega C(L)$ of the Cartan-Chevalley-Eilenberg
coalgebra $C(L)$. In fact, for a DG Lie algebra $L$ there exists 
a surjective quasi-isomorphism of DG algebras $\Omega C(L)\to U(L)$ 
(and even of DG Hopf algebras). Of course, 
$\Omega C(L)$ is much larger than $U(L)$: on the level of 
vector spaces the former
is isomorphic to tensor algebra $T^* \Lambda^*(L)$ on the exterior coalgebra   
$\Lambda^*(L)$,
while the latter is isomorphic to the symmetric coalgebra
$Sym^*(L)$ by PBW theorem. 

The DG algebra $\Omega C(L)$ also makes sense for a general $L_\infty$-algebra
$L$ and works well enough as a universal enveloping if we deal with 
DG algebras up to quasi-isomorphism. In other situations, one would like to 
have some structure on  $Sym^*(L)$ generalizing the usual universal 
enveloping. Since $A_\infty$-algebras relate to associative algebras 
 as $L_\infty$-algebras to Lie algebras, it is natural
to expect that $Sym^*(L)$ should be an $A_\infty$-algebra.

To construct it, we first consider a general $L_\infty$-algebra  $L$ as a 
DG vector space 
(= DG Lie algebra with trivial bracket). Then $C(L)$ turns into the 
symmetric coalgebra $Sym^*_c(sL)$ on the suspension $sL$ (isomorphic
as a vector
space to  the exterior coalgebra $\Lambda^*(L)$)
and the universal enveloping turns into the symmetric algebra 
$Sym^*_a(L)$.  Passing  from 
$\Omega Sym^*_c(sL)$ to $\Omega C(L)$ amounts to perturbing the 
differential on the tensor algebra and the standard techniques of
homological perturbation theory, cf. e.g. \cite{GLS}, 
give an $A_\infty$-structure on
$Sym^*(L)$. After the first draft of the present paper has been 
completed, it was brought to the author's attention that a similar
strategy (using filtrations instead of  perturbation
theory) was used by Polishchuk and Positselski in \cite{PP} in their
 proof of the PBW theorem.

However, functorial properties of such $A_\infty$-structure
will depend on a homotopy contracting $\Omega Sym^*_c(sL)$ onto 
$Sym_a^*(L)$.
For example, when $L$ is a finite dimensional vector space in degree
zero, one needs the homotopy to be $GL(L)$-invariant. 

This motivates a closer study of $\Omega Sym^*_c(sV)$ for a DG vector space
$V$. In Section 3 we prove an isomorphism of complexes, cf. Theorem 1:
$$
\Omega Sym^*_c(sV) \simeq k \oplus \bigoplus_{n \geq 1} \Big(V^{\otimes n}
\otimes_{k[\Sigma_n]} C_*(P_n)\Big)
$$
where $\Sigma_n$ is the symmetric group and $C_*(P_n)$ is the complex
computing the cell homology of the $n$-th permutahedron $P_n$
(the convex hull of the orbit of $(1, 2, \ldots, n)\in \mathbb{R}^n$
under the natural action of the symmetric group). 
 This leads to a
functorial - but not quite canonical - choice of a contracting homotopy. 

In Section 2 we construct the universal enveloping $U(L)$ and prove that
it has expected properties that generalize those of the classical 
 universal enveloping.
In particular, in Theorem 2 we show that $U(L)$ is a sort of 
``homotopy Hopf algebra"
even though the  operadic meaning of our construction
 remains unclear at the moment, see Section 4
for a  discussion. At the moment, $L \mapsto U(L)$ falls short of
being a functor: we are only able to prove that $U(\psi) \circ U(\phi)
= U(\psi \circ \phi)$ if one of the $L_\infty$-morphisms $\psi, \phi$ is
\textit{strict}. In Theorem 3 we generalize the classical 
complex $C(L)\otimes U(L)$ and prove a derived equivalence between 
$C(L)$ and $U(L)$ (one of the versions of the BGG correspondence). In
Theorem 4 we show that appropriate categories of $A_\infty$-modules
over $U(L)$ and $L_\infty$-modules over $L$, are equivalent. 
At the end of Section 2 we also discuss
the example which has been the original motivation for this paper: namely, 
if $S(X)$ is the homogeneous coordinate algebra of a toric complete 
intersection $X$, then its ``Koszul dual" $E(X)$ is precisely the universal 
enveloping algebra of the $L_\infty$-algebra $L$ which is defined using
the equations of $X$. 

In Section 4 we discuss some further  questions related to 
the Hopf property of $U(L)$, operads, etc.
The appendix contains
standard results on differential homological algebra and homological perturbation 
theory.

\medskip
\noindent
\textbf{Acknowledgements.} The author thanks Jim Stasheff and Pavlo 
Pylyavskyy for the useful discussions. This work was supported by 
the Sloan Research Fellowship.

\section{The Universal Enveloping}

\subsection{Signs and suspensions}

We consider complexes of vector spaces $k$ over a field of characteristic
zero. We use cohomological grading, to be denoted by superscripts,
in which differentials have degree +1. If
$V$ is a complex, its suspension $sV$ is
defined by $(sV)^p = V^{p+1}$, $d(sv) = - s(dv)$. 
In particular $\deg(sv) = \deg v - 1$. All tensor products
are over $k$ unless indicated otherwise. 
Throughout this paper we use the Koszul sign rule
$$
(F \otimes G)(a \otimes b) = (-1)^{\deg G\cdot \deg a}F(a) \otimes G(b)
$$
If $V$ is a graded vector space 
$Sym^*(V) = \oplus_{k \geq 0}Sym^k(V)$ will 
stand for its graded symmetric tensors, i.e.
$Sym^k(V)$ is the space of vectors in $V^{\otimes k}$ which 
are invariant with respect to the \textit{graded} action 
of the symmetric group $S_k$ (i.e. whenever two odd elements are 
permuted this leads to a change of sign). If we disregard the grading
and assume that $V$ has only even vectors (resp. only odd vectors)
this will become the usual space of symmetric (resp. antisymmetric)
tensors. Note that $Sym^*(V)$ has standard structures of a 
commutative algebra $Sym^*_a(V)$ and a cocommutative coalgebra
$Sym^*_c(V)$.

\subsection{Universal enveloping: case of Lie algebras and the general plan.}

Let $L$ be a DG Lie algebra. One way - perhaps a little exotic - 
to construct its universal enveloping algebra is outlined below.
We use the notions of the reduced bar construction $B(A)$
of an augmented DG  algebra $A$ (and its $A_\infty$-version),
reduced cobar construction $\Omega(C)$ of a coaugmented 
coalgebra $C$, and the  Cartan-Chevalley-Eilenberg coalgebra
$C(L)$ of a DG Lie algebra $L$ (and its $L_\infty$-version). 
All these are reviewed in the appendix, see also \cite{MSS}, 
\cite{K} and \cite{LM}. 
Also, $U(L)$ will stand for the classical universal enveloping 
of a DG Lie algebra $L$, and later its $L_\infty$-version.

Consider the natural projection $s^{-1}C(L) = s^{-1}Sym^*(sL) \to L$
and extend it to a morphism of 
algebras $\Omega C(L)  \to U(L)$. Direct computation shows 
that this is actually a morphism of DG bialgebras. By Theorem 22.9
and the first equality on page 290 in \cite{FHT}, it is also a quasi-isomorphism. In Section 3 we essentially re-prove this 
assertion.

We can turn this property inside out and use as a definition. 
First, consider  $L$ with 
the same differential but trivial Lie bracket. The above
construction gives a quasi-isomorphism
of DG algebras $\Omega C(L)  \to Sym^*_a(L)$.  
Bringing back the original bracket on $L$ will deform  the differential
on $C(L)$, and therefore the differential on $\Omega C(L)$. The general
machinery of perturbation theory, see \cite{GLS} and the 
 next subsection, gives a new DG algebra structure on 
$Sym^*(L)$ and a multiplicative projection from $\Omega C(L)$ 
onto $Sym^*(L)$ which is still a quasi-isomorphism. 
In Theorem 2 (v) we prove   that the new
structure on $Sym^*(L)$ is  precisely the universal enveloping $U(L)$ 
(identified by PBW theorem with $Sym^*(L)$ as a  \textit{coalgebra}). 

This approach also gives a recipe for a general $L_\infty$-algebra $L$, 
since an $L_\infty$-structure also gives a perturbation of the 
differential on $C(L)$ and we can carry out a similar procedure of
adjusting the product on $Sym^*(L)$. By \textit{loc. cit.} such 
adjustment 
in general leads to an $A_\infty$-structure on $Sym^*(L)$. 
As the  procedure depends on a choice of homotopy on 
$\Omega C(L)$ our construction will be based on the following result.

\begin{theorem} For a complex $V$ 
set $A(V) =\Omega Sym^*_c(sV)$, $E(V) = Sym^*_a (V)$. Let $f_V: A(V) \to E(V)$
be the multiplicative extension of the projection $s^{-1} Sym^{\geq 1}
(sV) \to V$, and $g_V: E(V) \to A(V)$ the map given by
composition of natural 
embeddings
$$
Sym^n(V) \hookrightarrow V^{\otimes n} \hookrightarrow T^*(V) \hookrightarrow T^*(s^{-1} Sym^*_c(sV)) = \Omega Sym^*_c(sV)
$$ 
Then $f_V g_V = 1$ and there exists a contracting homotopy $h_V: A(V)\to A(V)$ which satisfies
$$
1 - g_V f_V = d h_V + h_V d; \qquad
f_V h_V= 0; \qquad h_V g_V= 0; \qquad h_V h_V = 0
$$
and is functorial in the following sense: for every morphism of complexes
$\phi: V \to W$ the natural induced map $A(V) \to A(W)$ fits
into commutative diagram
$$
\begin{CD}
A(V) & @>>> & A(W) \\  
@V{h_V}VV && @V{h_W}VV \\
A(V) & @>>> & A(W)
\end{CD}
$$
Moreover, one can choose $h_V$ to commute with the algebra 
anti-involution $\iota_\Omega$
on $\Omega Sym^*_c(sV)$ which acts by $(-1)$ on the space of 
generators $s^{-1}Sym^*_c(sV)$.  
\end{theorem} 
The proof of this theorem is given in Section 3. We will see that 
such a homotopy $h_V$ (or rather a system of homotopies $V \mapsto h_V$)
is not unique but its choice depends of purely combinatorial data that has
nothing to do with $V$. 

\subsection{Universal enveloping: construction and first properties.}

\textit{Construction.} 

\bigskip
\noindent
Let $(L, \{l_i\}_{i \geq 1})$ be an $L_\infty$-algebra.
First consider the complex $(L, l_1)$ 
and set in Theorem 1 $V = L$ which gives
a contraction $(f_L, g_L, h_L)$ from $A(L) = \Omega Sym_c^*(sL)$ to 
$E(L)= Sym_a^* (L)$. From this we produce
a contraction of the free tensor coalgebra $T^*_c (s\overline{A(L)})$ 
onto  the free tensor coalgebra $T^*_c (s\overline{E(L)})$ (here and 
below $\overline{(\cdot)}$ denotes the augmentation ideal). 
Recall  coproduct 
on  $T^*_c (s\overline{A(L)})$:
$$
\Delta_B[a_1, \ldots, a_n] = 1 \boxtimes [a_1, \ldots, a_n] 
+ [a_1, \ldots, a_n] \boxtimes 1 + \sum_{i = 1}^{n-1} [a_1, \ldots, a_i]
\boxtimes [a_{i+1}, \ldots, a_n],
$$ 
and similarly for $T^*_c (s\overline{E(L)})$.
On $s\overline{A(L)}\subset T^*_c(s\overline{A(L)})$ we set 
$f'_L = s f_L s^{-1}$, $g'_L = s  g_L 
 s^{-1}$, $h'_L = - s h_L s^{-1}$, and define
a contraction on $(s\overline{A(L)})^{\otimes n} \subset 
T_c^*(s\overline{A(L)})$ by 
$$
F^\circ_L = (f'_L)^{\otimes n}, \quad G^\circ_L = (g'_L)^{\otimes n},
\quad 
H^\circ_L = \sum_{t = 1}^n (g'_L f'_L)^{\otimes (t-1)} \otimes 
h'_L \otimes 1^{\otimes (n-t)}.
$$
Observe that 
 $H^\circ_L$   
satisfies the \textit{coalgebra homotopy} condition
$$
\Delta_B H^\circ_L = (H^\circ_L\otimes 1 + G^\circ_L F^\circ_L 
\otimes H^\circ_L)\Delta_B.
$$
Denote by $\delta_L^\circ$ and $d_L^\circ$ the differentials of the 
two tensor coalgebras, 
respectively. 
By definition $BA(L)$ differs from $T^*_c(s\overline{A(L)})$ only in its 
differential, given by 
$$
\delta_L= \delta_L^\circ + t_{\mu} + t_{L}
$$
where $t_{\mu}$ is the part that encodes  the product
on the tensor algebra $A(sL)$ and $t_{L}$ is the perturbation which encodes
the  $L_\infty$-brackets $l_i, i \geq 2$ on $L$, cf Section 5.1 in the 
appendix. 
Using Proposition 7 in the appendix we obtain  a new
contracting homotopy 
$$
F_L = (F^\circ_L)_{t_{\mu} + t_L}, G_L= (G^\circ_L)_{t_{\mu} + t_L}, 
H_L= (H^\circ_L)_{t_{\mu} + t_L}
$$
from $B\Omega C(L)$ to $T^*_c(s\overline{E(L)})$ with its new
differential $d_L= (d_L^\circ)_{t_{\mu} + t_L}$.
Since $t = t_{\mu} + t_L$ is a coalgebra perturbation, i.e.
$$
\Delta_B t = (t \otimes 1 + 1 \otimes t) \Delta_B,
$$
the new differential $d_L$ is again a coderivation, 
$F_{L}, G_{L}$ are morphisms of DG coalgebras and $H_{L}$ is a coalgebra homotopy, cf. 
\cite{GLS}. 

\bigskip
\noindent
\textbf{Definition.} 
\begin{enumerate}

\item Denote by $U(L)$
the vector space $E(L) = Sym^*(L)$ with the $A_\infty$-structure 
$\{m_i\}_{i \geq 2}$ given by 
the above coalgebra differential $d_L$ on $T^*_c(s\overline{E(L)})$.
Then $(T^*_c (s\overline{E(L)}), d_L)$ tautologically turns 
into the cobar construction $BU(L)$ of $U(L)$.

\item If $L, M$ are two $L_\infty$ algebras and $\phi: C(L)\to C(M)$
is an $L_\infty$ morphism, cf. \cite{LM}, let $U(\phi) = F_M \;
B\Omega(\phi) \; G_L:
BU(L)\to BU(M)$.

\item If $\phi: C(L)\to C(M)$ and $\psi: C(M)\to C(N)$ are two 
$L_\infty$-morphisms, set $H(\phi, \psi) = F_N \; B\Omega(\psi)
\; H_M \; B \Omega(\phi) \; G_L: BU(L)\to BU(N)$. 
\end{enumerate}

\begin{theorem} Let $\phi: C(L) \to C(M)$ be an $L_\infty$-morphism
of $L_\infty$-algebras $L, M$ and $\phi_1: L\to M$ be its first
component. Then

\begin{enumerate}
\item $U(\phi)$ is an $A_\infty$-morphism from $U(L)$ to $U(M)$
and its first component $U(\phi)_1: U(L) = Sym^*(L)  \to Sym^*(M) = U(M)$ 
is given by  symmetrization of $\phi_1$. 

\item If $\phi: L \to M$ is a strict morphism of $L_\infty$-algebras, 
i.e. $\phi_i = 0$ for $i \geq 2$, then the same holds for 
$U(\phi)$, i.e. $U(\phi)_i = 0$ for $i \geq 2$. 

\item The standard coproduct 
$\Delta: Sym^*(L) \to Sym^*(L)\otimes Sym^*(L)$ is a strict morphism of $A_\infty$-algebras, if the latter is given an $A_\infty$-structure
via the natural isomorphism
$$
Sym^*(L) \otimes Sym^*(L) \simeq Sym^*(L\oplus L).
$$

\item If $\phi: C(L)\to C(M)$ and $\psi: C(M)\to C(N)$ are 
two $L_\infty$-morphisms then $$U(\psi \circ \phi) - U(\psi)
\circ  U(\phi) = d_{U(N)} H(\phi, \psi) + H(\phi, \psi) d_{U(L)}:
BU(L) \to BU(N).$$
Moreover, if at least one of the morphisms $\phi, \psi$ is strict, then
$H(\phi, \psi)= 0$.

\item Suppose that the 2-truncation $(L, l_1, l_2)$ is a DG Lie algebra.
Then $(U(L), m_1, m_2)$ is a DG algebra isomorphic to the usual universal 
enveloping of $(L, l_1, l_2)$.

\item Let $\iota: U(L) \to U(L)$ be the linear involution that
corresponds to the action of $(-1)^k$ on $Sym^k(L)$. Then 
$$
m_n \circ \iota^{\otimes n} 
=  \iota \circ m_n \circ \omega_n
$$
where $\omega_n$ is the permutation $\{1, \ldots, n\} 
\to \{n, \ldots, 1\}$. In other words,  $\iota$ is a strict 
morphism $U(L) \to U(L)^{op}$, where
$(\cdot)^{op}$ is the opposite $A_\infty$-structure.

\item Let $n \geq 2$ and $v_1, \ldots, v_n \in L\subset U(L)$. 
Let $Alt(v_1 \otimes \ldots \otimes v_n)$ be the graded 
antisymmetrization 
of $v_1 \otimes \ldots \otimes v_n$. Then 
$$
m_n(Alt(v_1 \otimes \ldots \otimes v_n)) = l_n (v_1, \ldots, v_n).
$$

\end{enumerate}
\end{theorem}
\textit{Proof of (i) - (iv).}
To prove (i) first observe that $F_M$ and $G_L$ are DG coalgebra
morphisms by \cite{GLS} and $B\Omega(\phi)$ is a DG coalgebra
morphism since $\phi$ itself is a DG coalgebra morphism. Therefore
$U(\phi): BU(L)\to BU(M)$ is a DG coalgebra morphism encoding an
$A_\infty$-morphism $U(L)\to U(M)$. To compute the first component we need
to evaluate $U(\phi)$ on $v \in \overline{U(L)} \subset BU(L)$. 
But, for such an element, all terms in $F_M$, $G_L$ which involve
perturbation of the differentials on $B \Omega C(L)$, $B\Omega C(M)$,
are identically zero, therefore $U(\phi)(v) = F_M^\circ B \Omega (\phi) 
G_L^\circ (v)$ and the latter map is precisely given by the symmetrization
$Sym(\phi)$  of $\phi$.

To prove (ii) we observe that for a strict morphism $\phi$ one 
has $H_M^\circ B \Omega (\phi) = B\Omega(\phi) H_L^\circ$ by Theorem 1. 
Using the explicit formulas of the Basic Perturbation Lemma,
$$
F_M= F^\circ_M (1 - X_M H_M^\circ); 
\quad G_L = (1 - H^\circ_L X_L) G^\circ_L; \quad
H_M = H^\circ_M (1 - X_M H_M^\circ)
$$ 
and the 
side conditions $H_M^\circ H_M^\circ = 0, F_M^\circ H_M^\circ = 0, H^\circ_L
G^\circ_L = 0$ we obtain
$$
F_M \circ B\Omega (\phi) \circ G_L
= F_M^\circ  \circ B \Omega (\phi) \circ G_L^\circ = B Sym(\phi)
$$
Part (iii) is an immediate application of (ii) to the diagonal map
$L\to L\oplus L$, $x \mapsto x \oplus x$ which is a strict morphism of 
$L_\infty$-algebras. 

Finally, the left hand side in part (iv) by definition is equal to 
$$
 F_N B\Omega(\psi) (1 - G_M F_M) B \Omega (\phi) G_L
= F_N B \Omega (\psi) (\delta_N H_M+ H_M \delta_N) B \Omega (\phi) G_L
$$
and the assertion 
follows since $F_N$, $B \Omega(\psi)$, $B \Omega (\psi)$ and $G_L$ are
morphisms of complexes. To prove the vanishing we observe that, by 
Theorem 1, $H_M^\circ B \Omega(\phi) = B\Omega(\phi) H^\circ_L$ if
$\phi$ is strict, and similarly for $\psi$. Now the side 
conditions and the formulas for $F, G, H$ finish the proof.

\bigskip
\noindent
\textit{Proof of (v).} First we assume that $L$ is a Lie 
algebra, i.e. all $l_i$ vanish for $i \geq 3$. The $A_\infty$-structure
on $E(L)= Sym^*(L)$ is given by the following differential on $T^*_c(sE(L))$:
$$
d_L
= d^\circ_L + F^\circ_L \Big(\sum_{i \geq 0} (-1)^i \big((t_{\mu}+ t_L) 
H^\circ_L \big)^i 
\Big)(t_{\mu} + t_L) G^\circ_L
$$
To simplify this expression we first introduce a ``geometric grading"
on $\Omega Sym^*_c(sL)$ by declaring that elements 
of $s^{-1} Sym^k(sL)$ have degree $(k-1)$, and extending to 
$\Omega Sym^*_c(sL)$ multiplicatively (we can agree that 
 $k \subset \Omega Sym^*_c(sL)$ has degree $(-1)$ but that will 
not be used in the proof). From the point of view Lemma 5 in Section 3, 
this grading corresponds
to dimension of the cells of permutahedra. We extend it to 
$B \Omega Sym^*_c(sL)$ in the obvious way (again, setting to $(-1)$
on the constants).

Then $t_L$ vanishes on elements of geometric degree 0 since those 
elements are 
products of linear symmetric tensors, and the bracket $l_2$ encoded
by $t_L$ needs two inputs. Since the image of $G^\circ_L$ belongs to 
the degree 0 part we will have $t_L G^\circ_L = 0$. Also, 
the proof of Theorem 1, cf. Section 3.2, implies that 
 $H^\circ_L$ increases
the geometric degree by 1, $t_L$ decreases by 1, 
$t_{\mu}$ preserves it, while $F^\circ_L$ vanishes on elements
of positive degree. Consequently, the above formula for the deformed 
differential simplifies to 
$$
d_L
= d^\circ_L+ F^\circ_L \Big(\sum_{i \geq 0} (-1)^i 
\big(t_L H^\circ_L \big)^i \Big) t_{\mu} G^\circ_L
$$
Since the terms responsible for a multiple product
$m_n: U(L)^{\otimes n} \to U(L)$ are those which contain $t_{\mu}$ 
exactly $(n-1)$ times, we see that the differential
on $U(L)$ is the same on $Sym^*(L)$ and all $m_n$ with $n \geq 3$ vanish.
Therefore $U(L)$ is a DG algebra. Denoting the usual symmetric product 
in $E(L)= Sym^*(L)$ by $*$ we also see that for 
$x, y \in Sym^*(L)$ homogeneous  in the geometric grading:
$$
m_2(x, y) = x * y + (\textrm{terms\ of\  lower\ geometric\  degree}).
$$
Therefore,
the subspace $L \subset U(L)$ generates $U(L)$ as an algebra. For $v, u
\in L$ an explicit computation shows
$$
m_2(v, u) = v * u + \frac{1}{2} l_2(v, u). 
$$
Denote for a moment by $U^{cl}(L)$ the classical universal enveloping.
The last formula gives a surjective DG algebra morphism $U^{cl}(L) \to U(L)$
which is easily seen to be an isomorphism by an inductive argument involving
natural filtrations on both algebras.

Next, we assume that the higher products $l_i$, $i\geq 3$ of $L$ are not 
necessariy zero. Then the pertrubation $\delta_L = 
\delta_L^\circ + t_\mu + t_L$ can be 
split as $(\delta_L^\circ + t_\mu + t^2_L) + (t_L - t^{(2)}_L)$ 
where $t^2_L$ is the term
coming from the bracket $l_2$. The expression in the first parenthesis 
has square zero since by assumption $(L, l_1, l_2)$ is a DG Lie algebra. 
Denote by $F', G', H'$ and $d_L'$ the data corresponding to the 
perturbation $\delta^\circ_L + t_\mu + t^{(2)}_L$ and set 
$t^{\geq 3}_L= t_L - t^2_L$.
By Proposition 7 in the appendix the $A_\infty$-structure of $U(L)$ 
corresponds to 
the perturbation of $F', G', H'$ and $d_L'$ by $t^{\geq 3}_L$. 
In particular,
the differential of $BU(L)$ is given by 
$$
d'_L + F'\big(\sum_{i \geq 0} (-1)^i (t^{\geq 3}_L H')^i \big) t^{\geq 3}_LG'.
$$
Evaluating the second term on $s\overline{U(L)} \subset BU(L)$ and 
$s\overline{U(L)}\otimes s \overline{U(L)} \subset BU(L)$ will give
zero for the following reasons. Firstly, for $x \in s\overline{U(L)}$ 
we have $G' (x) = G^\circ_L(x)$ since $(t_{\mu} + t^2_L)G^\circ_L(x) 
= 0$. But then $t^{\geq 3}_L G'(x) = t^{\geq 3}_L G^\circ_L(x) = 0$ since
$t^{\geq 3}$ vanishes on terms of geometric degree $\leq 1$. 
Secondly, for $x_1, x_2 \in s\overline{U(L)}$ by a similar computation
$$
G'(x_1 \otimes x_2) = \Big[\sum_{i\geq 0} (-1)^i(H^\circ_L t^2_L)^i\Big]
H^\circ_L (G^\circ_L(x_1) \otimes G^\circ_L (x_2))
$$
Since 
$H$ increases the geometric degree by 1 and $t^2_L$ 
decreases it 1, the above expression has geometric
degree 1, so $t^{\geq 3}_L$ vanishes on it.
This means that the differential and the product of $U(L)$ are the 
same as for the 2-truncation $(L, l_1, l_2)$, which finishes
the proof of (iv).

To prove (vi) for $n \geq 3$ consider a similar anti-involution $\iota_\Omega:
\Omega C (L) \to \Omega C(L)^{op}$ of Theorem 1. Let 
$\widehat{\omega}$ be a linear involution on $BU(L)$ which acts
by $\omega_n$ on $(s\overline{U(L)})^{\otimes n}$ and 
use the same notation for the corresponding involution on $B\Omega C(L)$.
Denote by 
$\pi: BU(L)\to \overline{U(L)}$ projection onto the first 
component.  Also, let $B\iota$, $B\iota_\Omega$ be
the linear involutions on the bar constructions
which act by $s^{\otimes n} \iota^{\otimes n} (s^{\otimes n})^{-1}$, 
$s^{\otimes n} \iota_\Omega^{\otimes n} (s^{\otimes n})^{-1}$
on the $n$-th tensor components, respectively. 
Since $\omega_n s^{\otimes n} = (-1)^{\frac{n(n-1)}{2}} 
s^{\otimes n} \omega_n: (\overline{U(L)})^{\otimes n} 
\to (s\overline{U(L)})^{\otimes n}$, we need to show that 
$$
(-1)^{\frac{n(n-1)}{2}} 
\pi (F^\circ_L X_L G^\circ_L) (B\iota\; \widehat{\omega})
= (B\iota\; \widehat{\omega}) \pi (F^\circ_L X_L G^\circ_L)
$$
on $(s\overline{U(L)})^{\otimes n}$. By Section 5.2 
in the appendix $X_L$ is a sum of several terms  of the form 
$$
(-1)^s a_1 \ldots a_s t_\mu
$$
where each $a_i$ is either $(t_L H^\circ_L)$ or $(t_\mu H^\circ_L)$.
If such a term is to give a nonzero contribution to the expression
above, the operator $t_\mu$ should be used exactly $(n-1)$ times,
since we need to get from $(s\overline{U(L)})^{\otimes n}$
to $s\overline{U(L)}$. It is easy to see that 
$$
\quad 
(B\iota\; \widehat{\omega})F^\circ_L = F^\circ_L (B\iota_\Omega\; \widehat{\omega});
\quad 
(B\iota_\Omega\; \widehat{\omega}) G^\circ_L = G^\circ_L (B\iota\; \widehat{\omega})
$$
and that $(B\iota_\Omega\; \widehat{\omega})$ commutes
with the operators $t_L$ and $H^\circ_L$.  Now what we need to 
prove follows from the following formula, easily checked by direct
computation:
$$
(B\iota_\Omega\;\widehat{\omega}) t_\mu = (-1)^{i-1} 
t_\mu (B\iota_\Omega\; \widehat{\omega}):
(s\overline{U(L)})^{\otimes i} \to (s\overline{U(L)})^{\otimes (i-1)}.
$$
For $n =2$ the same argument works for $(m_2 - *)$ where $*$ is the usual
product on $Sym^*(L)$. Since $*$ is commutative, the assertion holds for
$m_2$ as well. For $n=1$, the differential on $U(L)$ is the same as 
on $Sym^*(L)$ and the statement holds again.

\medskip
\noindent
Finally, (vii) is a restatement of Theorem 3 (i) below and 
its proof will be given there. $\square$

\subsection{Universal enveloping: categories of modules.}
Recall that $U(L)$ denotes the vector space $Sym^*(L)$ with 
the $A_\infty$-structure constructed in the previous subsection.
The next theorem deals with the notion of a generalized 
twisted cochain and the functors  defined by it, see appendix.
Part (iii) asserts a BGG-type equivalence to two derived
categories, $\mathcal{D} U(L)$ and $\mathcal{D}C(L)$.  
The derived category $\mathcal{D}U(L)$ is obtained by localizing
the category $Mod_\infty(U(L))$ of
strictly unital $A_\infty$-modules over $U(L)$ and
strictly unital morphisms (= the full subcategory of 
DG-comodules over $BU(L)$ which are free as comodules), at
the class of quasi-isomorphisms. The derived category $\mathcal{D}C(L)$
is obtained by localizing the category $Comodc(C(L))$ 
of cocomplete counital DG-comodules over $C(L)$, by the class
of weak equivalences (i.e. morphisms which induce a
quasi-isomorphism on the bar construction). See Chapter 2 in \cite{LH}
and Section 3.2 in \cite{B2} for more details.

\begin{theorem} The universal enveloping $U(L)$ has the following properties:
\begin{enumerate}
\item the composition $\tau: C(L)\to L \to U(L)$ is a generalized
twisted cochain;
\item the complex $C(L) \otimes_\tau U(L)$ is quasi-isomorphic to $k$ 
and the DG algebra morphism 
$\Omega C(L) \to \Omega B U(L)$ induced by $\tau$, is a 
quasi-isomorphism;
\item the functors $M \mapsto M\otimes_\tau C(L)$ and $N \mapsto N\otimes_\tau U(L)$ induce mutually inverse equivalences of the derived categories 
$\mathcal{D}C(L)$ and $\mathcal{D}U(L)$.
\end{enumerate}
\end{theorem}
\textit{Proof.} To prove (i), start with the composition 
$$
C(L)\to B\Omega C(L)\stackrel{F_L}\longrightarrow BU(L).
$$
Since it is a DG coalgebra morphism, by 5.3 in the appendix, its projection 
onto  $U(L)$ is a generalized twisted cochain $C(L)\to U(L)$. It is easy to 
check that it coincides with  $\tau$.

Part (ii) is known when $L$ is an abelian and the general case follows
by perturbation lemma as in the construction before Theorem 2. 
Alternatively, for the fist assertion 
we could first replace $U(L)$ by $\Omega C(L)$ 
where the corresponding results are again well known, cf. \cite{FHT},
and then pass from $\Omega C(L)$ to $U(L)$ using the strategy of 
\cite{AAFR}; while the second assertion is entirely similar to the
case of Lemma 6 in \cite{B2}. 

Part (iii) is a standard consequence of (ii), see Section 3.3.
of \cite{B2} and \cite{LH} for the associative case. $\square$

\bigskip\noindent
We can also construct a pair of functors relating $L$-modules to 
$U(L)$-modules. 
Let $Mod(L)$ be the category of $L_\infty$-modules over $L$ and
$L_\infty$-morphisms (= the category of DG comodules over $C(L)$
which are free as $C(L)$-comodules).  
By the appendix, we can also view an $L$-module structure on $M$ 
as a twisted
cochain $\tau: C(L)\to End(M)$. The corresponding DG coalgebra map
$C(L)\to B End(M)$ admits a canonical factoring
$$
C(L)\to B\Omega C(L)\to B End(M)
$$
since we can extend $\tau$ to a DG algebra map $\Omega C(L) \to End(M)$
and then apply the bar construction. Therefore, composing with 
$G_L: BU(L)\to B \Omega C(L)$ we get a DG-coalgebra map 
$BU(L)\to B End(M)$, i.e. a strictly unital $A_\infty$-module structure
on $M$. This defines a functor 
$$
\mathcal{G}: Mod(L) \to Mod_\infty(U(L))
$$
In the opposite direction, we start with a DG coalgebra morphism
$BU(L)\to B End(M)$ and then composing with the canonical map
$C(L) \to B \Omega C(L)$ and $F_L: B\Omega C(L)\to BU(L)$ we get
a DG coalgebra map $C(L) \to B End(M)$, i.e. a twisted cochain 
$C(L)\to End(M)$ which gives $M$ a structure of 
an $L_\infty$-module over $L$. 
This defines a functor 
$$
\mathcal{F}: Mod_\infty(U(L)) \to Mod(L).
$$
Observe that in both cases the underlying vector space does not
change. 
\begin{theorem}
The above functors $\mathcal{G}$, $\mathcal{F}$ are mutually
inverse equivalences.
\end{theorem}
\textit{Proof.} In one direction, supppose we start with an 
$A_\infty$-module structure on $M$ given by  $BU(L) \to B End(M)$. 
Applying $\mathcal{G} 
\mathcal{F}$ amounts to considering the composition
$$
BU(L)\stackrel{G_L}\longrightarrow B\Omega C(L)
\stackrel{F_L} \longrightarrow BU(L) \to B End(M).
$$
Since the composition of the first two arrows is identity, we
conclude that the identity map on $M$ gives an isomorphism of
$A_\infty$-modules  $\mathcal{G} \mathcal{F}(M)$  and $M$. 

In the other direction, suppose we start with a twisted cochain
$C(L)\to End(M)$ and construct $B \Omega C(L)\to B End(M)$ as
above. The $L_\infty$-module corresponding to $\mathcal{F} 
\mathcal{G} (M)$ is obtained from a DG coalgebra morphism
$$
C(L)\to B\Omega C(L) \stackrel{F_L}\longrightarrow BU(L)
\stackrel{G_L}\longrightarrow B\Omega C(L) \to B End(M)
$$
In view of $G_L F_L= 1 - \delta_L H_L - H_L \delta_L$ it 
suffices to show that the composition 
$$
C(L)\to B\Omega C(L) \stackrel{\delta_L H_L + H_L \delta_L}
\longrightarrow  B\Omega C(L) \to B End(M)
$$ 
is zero. That in its turn would follow from the vanishing of 
$$
C(L)\to B\Omega C(L) \stackrel{H_L}\longrightarrow 
 B\Omega C(L).
$$
But the latter holds since $h_L$ vanishes on 
$s^{-1}\overline{C}(L) \subset \Omega C(L)$ by its
construction, see Section 3.2 (the homotopy $\mathcal{H}_n$ 
vanishes on the 
top-dimensional cell of the permutahedron $P_n$). 
Thus, the idenitity on $M$ also gives an isomorphism of 
$L_\infty$-modules $M$ and $\mathcal{F}\mathcal{G}(M)$, which 
finishes the proof.
$\square$

\subsection{An example: toric complete intersections.}

The following example had originally motivated our study of
$L_\infty$-algebras. See \cite{B1} and \cite{B2} for details. 

Let $X \subset \mathbb{P}^\Sigma$ be 
a complete intersection in a toric variety defined by a fan $\Sigma$.
Then $X$ has a ``homogeneous coordinate ring" $S(X) = Sym^*(V)/J$,
a quotient of a polynomial ring by an ideal generated by a regular
sequence of polynomials $W_1, \ldots, W_m$. For a general toric 
variety $S(X)$ will be graded by a finitely generated abelian group 
$A(X)$ and $W_1, \ldots, W_m$ will be homogeneous in this grading
(but not the usual grading of $Sym^*(V)$). One can always assume that
$W_1, \ldots, W_m$ have no linear terms.

In this setting, define the ``Koszul dual" of $S(X)$ as
the Yoneda algebra $E(X) = Ext^*_{S(X)} (k, k)$ with its natural 
$A_\infty$-structure (defined in general up to $A_\infty$-homotopy). 

Introducing formal degree 2 variables $z_1, \ldots, z_m$ which span 
a vector space $U$ we 
can define  an $L_\infty$-algebra 
$L = s^{-1} V^\vee \oplus U$ by viewing the formal sum $W= \sum W_i (sz_i)$
as a differential on $C(L) = Sym^*_c(V^\vee\oplus sU)$, if we agree
that $W_j$ act by differential operators on $Sym^*(V^\vee)$.

It was shown in \cite{B1} and \cite{B2} that the Koszul dual $E(X)$ may be 
identified with the universal enveloping $U(L)$ (the two papers 
quoted used Koszul type-resolutions instead of $\Omega C(L)$ which 
still lead to the same $A_\infty$-structure, perhaps after
a change of contracting homotopy). The interpretation 
in terms of $Ext$ groups also follows from Theorem 3 (ii). 

\section{A homotopy on the cobar construction}
\subsection{Permutahedra}

Let $n \geq 1$ and $1 \leq d \leq n$ and set $P(n, d)$ to be the
set of ordered partitions of $\{1, \ldots, n\}$ which have $d$
parts. Equivalently, any such partition can be viewed as a
surjective map $\psi: \{1, \ldots, n\} \to \{1, \ldots, d\}$: 
setting $\psi_i = \psi^{-1} (i)
\subset \{1, \ldots, n\}$, $1 \leq i \leq d$ we get an 
ordered partition $[\psi_1| \ldots | \psi_d]$.
There exists, cf. e.g. \cite{SU},
 a polytope $P_n$, called the $n$-th permutahedron, 
such that  $P(n, d)$ labels the faces of dimension $n-d$
in $P_n$. In particular, $P_n$ has dimension $n-1$ and its vertices
are labeled by permutations of $\{1, \ldots, n\}$. 

To consider the homology complex of $P_n$ 
define an orientation of 
$$
\psi: \{1, \ldots, n\} \to \{1, \ldots, d\} 
$$
as an equivalence class of orderings on  each subset $\psi_i$, such 
that two orderings are equivalent  if they differ by an even 
permutation of $\{1, \ldots, n\}$. We choose the  
orientation  corresponding to the 
natural increasing ordering on  $\psi_j$. 

Let $C_*(P_n)$ be the \textit{homology} complex of $P_n$
with grading inverted to ensure that differential has degree $+1$
(thus, $C_*(P_n)$ is concentrated in degrees $-n+1, \ldots, 0$). 
The notation $\psi = [\psi_1 | \ldots | \psi_d]$ allows to reduce 
most of the signs below to the Koszul sign rule if we assume 
that the symbol $|$ has degree $(+1)$ and each of the elements in 
$\psi_i$ degree $(-1)$. 

The differential of $C_*(P)$, cf. \cite{SU}, is given by:
$$
\partial [\psi_1| \ldots |\psi_d] 
= \sum_{\genfrac{}{}{0pt}{}{1 \leq k \leq d}{M \varsubsetneq\psi_k}}
(-1)^{\psi, M}
[\psi_1 | \ldots | \psi_{k-1} | M | \psi_k \setminus M |
\psi_{k+1} |\ldots | \psi_d].
$$
The sign is 
$$
(-1)^{\psi, M} = (-1)^{m_1 +\ldots + m_{k-1} + (k-1) + \#M} (-1)^{\sigma_M} 
$$
where $m_i = \#\psi_i$ and $\sigma_M$ is the unshuffle that takes
$\psi_k$ to $[M | \psi_k \setminus M]$ (again, taken  with the 
natural increasing ordering). 
The symmetric group $\Sigma_n$ acts from the left
on each $P(n, d)$ and on $C_*(P_n)$:
$$
\sigma [\psi_1| \ldots | \psi_d] = \pm 
[\sigma(\psi_1)| \ldots | \sigma (\psi_d)].
$$ 
where the sign is $(+1)$ if the ordering induced from $\psi$ by 
$\sigma$ is equivalent to the increasing ordering, and $(-1)$
otherwise. In addition, $C_*(P_n)$ has an involution
$$
\nu_n [\psi_1| \psi_2| \ldots | \psi_{d-1}| \psi_d] = -
(-1)^{n (d-1) +\frac{(d-1)(d-2)}{2} + \sum_{i <j} m_i m_j}
[\psi_d | \psi_{d-1}|  \ldots | \psi_2 |\psi_1]
$$
which commutes with the differential and the $\Sigma_n$-action.
Therefore, we actually have a  $\Sigma_n \times \mathbb{Z}_2$-action
on $C_*(P)$. 

\bigskip
\noindent
Define a bilinear map $\Theta: V^{\otimes n} \times C_*(P_n)
\to \Omega Sym^*_c(sV)$ by 
$$
\Theta(v_1 \otimes \ldots \otimes v_n, [\psi_1| \ldots | \psi_d]) = 
(-1)^{(n-d)
(\sum_i \deg v_i)}
s^{-1} (s^{\otimes m_1}) \otimes \ldots \otimes s^{-1}(s^{\otimes m_k})
\big[(v_1 \otimes \ldots \otimes v_n)\cdot \sigma_\psi\big]
$$
where $\sigma_\psi$ is the permutation $\{1, \ldots, n\}
\to [\psi_1 | \ldots | \psi_d]$ and each $s^{-1}(s^{\otimes m})$
is viewed as a map $V^{\otimes m} \to s^{-1} Sym^m(sV)$, 
$u_1 \otimes \ldots \otimes u_m \mapsto \pm s^{-1} (su_1 \ldots su_m)$
with the sign determined by the Koszul  rule.  The following lemma
amounts to a direct computation.

\begin{lemma}
The map $\Theta$ induces an
isomorphism of complexes
$$
\Omega Sym^*_c(sV)\simeq k \oplus \bigoplus_{n \geq 1} 
\Big(V^{\otimes n} \otimes_{k[\Sigma_n]} C_*(P_n)\Big)
$$
which takes $\iota_\Omega$ to $1 \oplus \bigoplus_{n \geq 1} (1 \otimes \nu_n)$.
\end{lemma}
\subsection{Proof of Theorem 1.}

Since $P_n$ is a convex polyhedron, the complex $C_*(P_n)$
has cohomology $k$ in degree 0, and zero everywhere else.
Let $\mathcal{F}_n: C_*(P_n) \to k$, $\mathcal{G}_n: k \to C_*(P_n)$,
be the natural $\Sigma_n \times \mathbb{Z}_2$-equivariant projection 
and embedding,
respectively (where $k$ is viewed as a trivial 
$\Sigma_n \times \mathbb{Z}_2$-module).
Since we are working in characteristic zero, we can find a 
$\Sigma_n \times \mathbb{Z}_2$-equivariant contracting homotopy 
 $\mathcal{H}_n: C_*(P_n)\to C_*(P_n)$. 
It is well known,  see e.g. Section 2.1 in \cite{LS}, 
that we can also assume  the side conditions:
$$
\mathcal{H}_n \mathcal{G}_n = 0, \quad \mathcal{F}_n \mathcal{H}_n = 0, 
\quad \mathcal{H}_n \mathcal{H}_n = 0
$$
(if the first two identities are not satisfied then replace $\mathcal{H}_n$
by $\mathcal{H}_n' = (1 - \mathcal{G}_n \mathcal{F}_n) \mathcal{H}_n 
(1 - \mathcal{G}_n \mathcal{F}_n)$, then if the last identity is not 
satisfied, replace $\mathcal{H}'_n$ by $\mathcal{H}''_n = 
\mathcal{H}'_n d \mathcal{H}'_n$; these explicit formulas also show 
that equivariance will still hold).

Using the decomposition of the previous lemma, set 
$$
h_V = 0 \oplus \bigoplus_{n \geq 1} (1 \otimes \mathcal{G}_n)
$$ 
By  the $\Sigma_n \times \mathbb{Z}_2$-equivariance
 it follows that $h_V$ is a homotopy contracting 
$\Omega Sym^*_c(sV)$ to 
$$
k \oplus \bigoplus_{n \geq 1} \Big(V^{\otimes n} \otimes_{k[\Sigma_n]} k\Big)
= Sym^*(V)
$$
and that $h_V$ commutes with the anti-involution $\iota_\Omega$ as
well. 

\subsection{Relation with semistandard tableaux.}

Our original approach to  Theorem 1 was based on the
equivalent language of semistandard tableaux. The main advantage
of using permutahedra is better compatibility with 
the involution $\iota_\Omega$ on $\Omega Sym^*_c(sV)$. On the 
other hand, semi-standard tableax give an explicit
decomposition of $\Omega Sym^*_c(V)$ into irreducible $GL(V)$-modules
(e.g. when $V$ is a finite dimensional vector space in homological
 degree 0). 
 These results (perhaps known to experts in combinatorics) are not 
used in this paper, and the proof
is left to the interested reader.

The link between permutahedra and Young tableaux becomes clear if
we consider the faces of  $P_n$ which correspond
to ordered partitions $\psi = [\psi_1 | \ldots | \psi_d]$ with 
fixed $m_i = \# \psi_i$. Denoting $\textbf{m} = (m_1, \ldots, m_d)$ we
see that the set of such faces is a single $\Sigma_n$-orbit 
of 
$$
\psi_{\textbf{m}} = [1, \ldots, m_1 | (m_1 + 1), \ldots, (m_2 + m_1)| \ldots | 
(m_1 + \ldots + m_{d-1} + 1), \ldots, n]
$$
If orientations are taken into account, it becomes clear that the line
$k \cdot \psi_{\textbf{m}} \subset C_*(P_n)$
is isomorphic to the sign representation $\rho_{\textbf{m}}$ of the 
stabilizer $\Sigma_\textbf{m} = \Sigma_{m_1} \times \ldots \times 
\Sigma_{m_d} \subset \Sigma_n$. Therefore, the $\Sigma_n$-submodule
$$
M_{\textbf{m}} = \bigoplus_{\{\psi| \#\psi_i = m_i \; \forall i\}} 
k \cdot \psi \subset C_*(P_n)
$$
is the induced representation 
$\rho \uparrow_{\Sigma_\textbf{m}}^{\Sigma_n}$. If $\mathcal{S}^\lambda$
is the irreducible  Specht module corresponding to a partition $\lambda$,
cf. e.g. \cite{S},
its multiplicity in $M_{\textbf{m}}$ can be computed as
the number of column-semistandard tableaux $T$ with content 
$\textbf{m}$, cf. Theorem 2.11.2 in \textit{loc. cit.} 
Thus, Lemma 5 above will give a decomposition of 
$\Omega Sym^*_c(sV)$ in terms of Schur complexes.  It takes 
some additional effort to make all explicit homomorphisms 
compatible with the differential.  

Let $\lambda$ be a partition of $n$ and use the same notation for 
the corresponding Young diagram. Choose a $\lambda$-tableau $T$, i.e.
a bijective map $\{\lambda \} \to \{1, \ldots, n\}$ where $\{\lambda\}$
is the set of cells in $\lambda$. Let $C_T, R_T \subset \Sigma_n$
be the column stabilizer and row stabilizer, respectively, i.e. those
permutations which preserve values in the columns, resp. rows of $T$. 
Setting
$$
c_T= \sum_{\sigma \in C_T}\sigma; \qquad r^-_T = \sum_{\sigma \in R_T} 
(-1)^\sigma \sigma; \qquad e_T = c_T r^-_T
$$ 
we can define the Schur complex $S^T(V) = (V^{\otimes n}) e_T$ for
any complex of vector spaces $V$. Here we use alternation in the
rows of $T$, rather than columns, because of the suspension $sV$
involved in $\Omega Sym^*_c(sV)$. 

Now suppose that $T$ is standard, i.e. the values increase in rows
and columns. Set 
$$
J_T = \{i\ | \ 1 \leq i \leq n-1, \textrm{\ and\ } T^{-1}(i)
\textrm{\ is\ strictly\ above\ }
T^{-1}(i+1)\} 
$$
By construction $J_T\subset \{1, \ldots, n-1\}$. For any $J\subset J_T$
with $p$ elements consider the unique weakly increasing surjective map
$$
\zeta_J:\{1, \ldots, n\} \to \{1, \ldots, n-p\} 
$$
such that $J= \{i\ |\  \zeta_{J} (i) = \zeta_{J}(i+1)\}$. Then
the composition 
$$
T_J: \{\lambda \} \stackrel{T}\longrightarrow \{1, \ldots, n\} 
\stackrel{\zeta_J}\longrightarrow\{1, \ldots, n-p\}
$$
is a \textit{column-semistandard tableaux}, i.e. the values 
increase weakly 
in the columns and strictly  in the rows.
It is easy to see that every surjective map $U: \{\lambda \} 
\to \{1, \ldots, n-p\}$ which is a
column-semistandard tableau, has the form
$T_J$ for unique $T$ and $J\subset J_T$.

\begin{theorem} One has a direct sum decomposition 
$$
\Omega Sym^*_c(sV) \simeq k \oplus \bigoplus_{\lambda} 
\bigoplus_{\genfrac{}{}{0pt}{}{T \textrm{\ is\ a\ standard}}{
\lambda-\textrm{tableau}}}
\big(C_T\otimes S^T(V)\big)
$$
where $C_T$ is a combinatorial complex spanned in degree $(-p)$ by 
$T_J$ with $J\subset J_T, \#J = p$ and differential given by 
$$
\partial (T_J)= \sum_{j \in J} (-1)^{\#X(J, j)}\; T_{(J \setminus j)}; 
\qquad X(J, j)= \{i\;|\; 1 \leq j \leq j-1, i \notin J\}
$$
\end{theorem}
To describe the isomorphism explicitly, for any $J \subset J_T$ let
$\textbf{m}(J) = (m(J)_1, \ldots, m(J)_{n-p})$ with
 $m(J)_i = \#T_J^{-1}(i)$ and $\sigma_{\textbf{m}(J)}
\in k [\Sigma_n]$ the average of all elements in the corresponding
subgroup $\Sigma_{\textbf{m}(J)} \subset \Sigma_n$. Then for 
$u \in S^T(V)
= V^{\otimes n} e_T$ we set 
$$
(T_J\otimes u) \mapsto
\frac{1}{m(J)_1! \ldots m(J)_{n-p}!}
\pi_J (u \sigma_{\textbf{m}(J_T)})\in \Omega Sym^*_c(sV)
$$ 
where $\pi_J$ is the composition 
$$
V^{\otimes n} \to Sym^{m(J)_1} (V) \otimes \ldots \otimes 
Sym^{m(J)_{n-p}}(V) 
\to s^{-1} Sym^{m(J)_1} (sV) \otimes \ldots \otimes 
s^{-1} Sym^{m(J)_{n-p}} (sV)
$$
Note that the complex $C_T$ may be indentified with the standard Koszul
complex on the vector space with basis labeled by elements of $J_T$;
and one could use this identification to write \textit{explicitly} a 
homotopy $h_V$
satisfying the functoriality condition of Theorem 1.
For example one could use the following homotopies on $C_T$:
$$
h_T(T_J) = \frac{1}{\#J_T} \sum_{j \in (J_T\setminus J)} (-1)^{\#X(J, j)}
T_{(J\cup j)}
$$
However, to ensure
that $h_V$ commutes with the involution $\iota_\Omega$ we may have
to replace it by $h'_V = \frac{1}{2}(h_V + \iota_\Omega h_V \iota_\Omega)$
and this has no apparent meaning in terms of semistandard tableaux. 

\section{Further questions}

\begin{itemize}

\item The present approach, in principle, should 
give a combinatorial formula for 
the product in the usual universal enveloping. It would be interesting
 to write it  explicitly. 

\item The homotopy $h_V$ constructed in Section 3 is not  
canonical: one is still making choices in terms of 
faces of the permutahedra. Is there a special choice of 
$h_V$ which results in any additional properties of $U(L)$, e.g.
$U(\psi) \circ U(\phi) = U(\psi \circ \phi)$ for all $\psi$, $\phi$? 

\item It would be nice to have a more thorough understanding of
the correspondence $L \mapsto U(L)$ from the operadic point of view. Also,
in Theorem 2 we prove that the diagonal map of $Sym_c^*(L)$ is 
a strict morphism of $A_\infty$-algebras if $Sym^*(L) \otimes Sym^*(L)$
is identified with $Sym^*(L \oplus L)$. Although the  ``$A_\infty$ tensor
product"
$$
Sym^*(L)\otimes Sym^*(M)\simeq Sym^*(L \oplus M)
$$
is extremely natural in the present
context, its relation to  such operadic constructions as the 
Saneblidze-Umble diagonal on permutahendra and associahedra, cf. 
\cite{SU}, or  the diagonal on the W-construction of the associative operad,
cf. \cite{MS}, remains a mystery for the author. 

The $W$-construction may be relevant,
since for any \textit{associative} algebra $A$  a contraction 
$(F, G, H)$ from $A$ to $E$ defines a $W$-algebra structure
on  $E$: using the terminology of \cite{MS}, non-metric edges
will be labeled by $GF$, metric edges by $H$ and internal vertices by 
multiple products. The usual $A_\infty$-structure, cf.  \cite{GLS},
\cite{KS}, is induced via the operadic map $A_\infty \to W$ described
in  \cite{MS}. Moreover, the diagonal on the $W$-construction,
cf. \textit{loc. cit.}, corresponds precisely to the tensor product of
homotopies $(H_1, H_2)\mapsto H_1 \otimes 1 + G_1 F_1 \otimes H_2$.

\item In a recent spectacular work, cf. \cite{M}, 
Merkulov has proved that a general 
homotopy  Lie bialgebra can be quantized, i.e. it defines 
 a homotopy bialgebra 
structure on $Sym^*(L)$. This construction involves some non-explicit
choices of operadic maps and a more transparent version of it is highly 
desirable. Is it possible to describe this quantization along the lines
of Kazhdan-Etingof using Theorem 4 in this paper?

\end{itemize}

\section{Appendix}

\subsection{Standard constructions of differential homological algebra.}

Let $L$ be a DG Lie algebra with differential 
$l_1$ and the bracket $l_2: L^{\otimes 2} \to L$.
Its \textit{Cartan-Chevalley-Eilenberg construction} $C(L)$ is the 
DG coalgebra $Sym^*_c(sL)$ with the differential 
$\delta_C = c_1 + c_2$ defined as follows. Let 
$s^{\otimes n}: L^{\otimes n} \to (sL)^{\otimes n}$ be the obvious 
degree $(-n)$ isomorphism and set 
$$
c_1 = - s\; l_1 s^{-1}:sL\to sL; \quad
c_2 = s\; l_2 (s^{\otimes 2})^{-1}:(sL)^{\otimes 2} \to sL
$$
extending these maps to $Sym^*_c(sL)$ as coderivations. 
The property $\delta_C^2 = 0$
follows from $l_1^2 = 0$, the Leibniz Rule and the Jacobi 
Identity for $l_2$. 

\bigskip
\noindent
In general, if $L$ is a graded
vector space and $\delta$ is a differential on $Sym^*_c(sL)$ 
which is a coderivation,
we can consider compositions $c_n: (sL)^{\otimes n} 
\to Sym^*_c(sL) \stackrel{\delta}\longrightarrow Sym_c^*(sL)\to sL$
and define $l_n: L^{\otimes n} \to L$ via 
$$
c_n = (-1)^n s\; l_n (s^{\otimes n})^{-1} 
$$
Then $\{l_n\}_{n\geq 1}$ give $L$ the structure of an 
$L_\infty$-algebra, cf. \cite{LM}. If $\phi: (Sym^*_c(sL), \delta)
\to (Sym^*_c(sL'), \delta')$ is a (degree zero) morphism of DG 
coalgebras, 
by a similar formula we get a sequence of degree $1-i$ 
maps $\phi_1:
\Lambda^i(L)\to L'$. The sequence $\{\phi_i\}_{i \geq 1}$ 
(or, equivalently, the original morphism $\phi$) is called
an $L_\infty$-morphism from $L$ to $L'$. 

The main purpose of this article is 
to provide a construction of the universal enveloping for $L$.
To that end, we need two more definitions.

\bigskip\noindent
Let $A = k \oplus \overline{A}$ be an augmented DG algebra with 
differential $m_1$ and product $m_2$. Its
\textit{reduced cobar construction} $B(A)$ is the tensor 
coalgebra $T^*_c(s\overline{A})$ with the differential
$\delta_B= b_1 + b_2$ defined in a similar way: 
$$
b_1 = (-1)^1 s\; m_1 s^{-1}: s\overline{A} 
\to s\overline{A}; \quad b_2 = (-1)^2 s \;m_2 (s^{\otimes 2})^{-1}:
(s\overline{A})^{\otimes 2} \to s\overline{A}.
$$
Then
$b_1$ and $b_2$ extend uniquely to $B(A)$ as coderivations and 
$\delta_B^2 =0$ follows from $m_1^2 =0$, the Leibniz Rule and 
associativity of $m_2$.

Again, one can consider a general differential $\delta_B$ on $T^*_c(s\overline{A})$ which is a coderivation, and obtain operations
$m_n: A^{\otimes n}\to A$ by first considering
$$
b_n: (s\overline{A})^{\otimes n} 
\to BA \stackrel{\delta_B}\longrightarrow BA \to s\overline{A}
$$
and then writing 
$$
b_n = (-1)^n s\; m_n (s^{\otimes n})^{-1}.
$$
The resulting operations $\{m_n\}_{n \geq 1}$
give $A$ a structure of an $A_\infty$-algebra, cf. \cite{K}, 
\cite{MSS}. Since we use the reduced bar construction, $A$ is 
automatically \textit{strictly unital}, i.e.
$$
m_n(v_1, \ldots, v_n) = 0; \textrm{\ if\ } n \geq 3
\textrm{\ and \ } v_i = 1\textrm{\  for\ some\ } i
$$
and $m_2(v, 1) = m_2(1, v) = v$. If $f: (BA, \delta_B) 
\to (BA', \delta_B')$ is a DG coalgebra morphism, we get 
a sequence of degree $(1-i)$ maps $f_i: A^{\otimes i} \to A'$ which 
we call  an $A_\infty$-morphism from $A$ to $A'$. Again, 
since we use reduced bar constructions, the \textit{morphism} is 
automatically \textit{strictly unital}: $f_i = 0$ if $i \geq 2$ 
and one of its arguments is equal to $1 \in A$. 

\bigskip
\noindent
Finally, let $C= k \oplus \overline{C}$ be a coaugmented DG coalgebra. 
Its \textit{reduced cobar construction} is a DG algebra $\Omega(C)
= T^*_a(s^{-1}\overline{C})$ with the differential $\delta_\Omega
= \omega_1 + \omega_2$ where $\omega_1$ and  $\omega_2$ are 
obtained from the differential on $C$ and the 
reduced coproduct $\overline{\Delta}: \overline{C}\otimes 
\overline{C} \to \overline{C}$, respectively, using the same pattern 
as before (except this time $w_1$ and $w_2$ are 
extended from $s^{-1}\overline{C}$ to 
$\Omega(C)$ as derivations). If $C$ is 
\textit{cocommutative} the DG algebra $\Omega(C)$ also 
has a shuffle  coproduct $\Delta_\Omega: \Omega(C)
\to \Omega(C)\boxtimes\Omega(C)$ defined
on $s^{-1}\overline{C} \subset \Omega(C)$ by 
$$
\Delta_{\Omega}(u)= u \boxtimes 1 + 1 \boxtimes u
$$
and extended to $\Omega(C)$ multiplicatively. Thus, $\Omega(C)$
becomes a DG bialgebra (the fact that $\delta_\Omega$ is 
also a coderivation uses cocommutativity of $C$).

\subsection{Perturbation Lemma}

We recall the main result behind the pertrubation 
machinery, cf. e.g. \cite{Br}.
Let $(M, d_M)$,  $(N, d_N)$ be two complexes. Consider a contraction
$$
F: N \to M; \qquad G: M\to N; \qquad H: N\to N
$$
which satisfies the usual identities $FG= 1_N$, $1_M - GF = d_NH + Hd_N$.
Adjusting the homotopy $H$ if necessary, cf. Section 2
in \cite{LS},  one can always assume that 
the following ``side conditions" are also safisfied
$$
F H= 0; \qquad H H = 0; \qquad HG= 0.
$$
Now suppose  we are given a new differential $d_N +t$ on $N$ such 
that  $(tH)$ is locally nilpotent (i.e. for any element
$n \in N$ there is a positive integer $k(n)$ such that 
$(tH)^{k(n)} (n) = 0$). Then the infinite sum 
$$
X = t - t H t + tHt Ht - \ldots
$$
is well-defined. Introduce
$$
F_t = F(1 - XH); \quad G_t = (1 - HX)G; \quad H_t = H - HXH; \quad (d_M)_t = d_M + FXG
$$
\begin{prop} (Basic Perturbation Lemmma)
Under the assumptions introduced
above, $(F_t, G_t, H_t)$ is a contraction of the complex
$(N, d_N + t)$ to the complex $(M, (d_M)_t)$ which also satisfies
the side conditions.
\end{prop}
The following result can be checked by direct computation.
\begin{prop}
Suppose that $d_N + t_1$ and $d_N+ t_1 + t_2$ are pertrubations 
satisfying the above nilpotency condition. Then $(d_M)_{t_1 + t_2} 
= ((d_M)_{t_1})_{t_2}$ and similarly for   $F, G$ and $H$.
\end{prop}

\subsection{(Generalized) twisted cochains}

Let $C$ be a coagumented DG coalgebra and $A$ an augmented DG algebra  
A degree +1 map $\tau: \overline{C} \to \overline{A}$ is 
called a \textit{twisted cochain} if $\tau$ satisfies
$$
\tau d_C + d_A \tau = \mu \circ (\tau \otimes \tau) \circ \Delta
$$
where $\Delta: \overline{C} \to \overline{C} \otimes 
\overline{C}$ is the \textit{reduced} coproduct of $C$ and $\mu$ is 
the product in $A$.
This conditions guarantees that $\tau$ both the canonical  
coalgebra morphism $C \to BA$ and the canonical algebra morphism
$\Omega C\to A$, which extend $\tau$, commute with differentials.

When $A$ is a strictly unital $A_\infty$-algebra, one can write
a \textit{generalized twisted cochain} condition for $\tau$, cf. 
\cite{AAFR}: 
$$
\tau d_C + d_A \tau = \sum_{i \geq 2} \mu_i \circ \tau^{\otimes i} \circ \Delta^{(i)}
$$
where $\mu_i$ are the products in $A$ and $\Delta^{(i)}: \overline{C}
\to \overline{C}^{\otimes i}$ is the 
iteration of the reduced coproduct. This condition is equivalent to 
requiring that $C\to BA$ is a morphism of complexes.

Finally, if $L$ is an $L_\infty$-algebra then an $L_\infty$-module
structure on a vector space $M$ is defined by choosing a differential
$d$ on $C(L)\otimes M$ which makes it a DG-comodule over $C(L)$.
This differential encodes maps $\Lambda^k(L) \otimes M \to M$ which 
satisfy a series of quadratic identities arising from $d^2 = 0$.
It follows from the definitions that the same structure is also encoded
by a twisted cochain $C(L)\to End(M)$. Similarly, $A_\infty$-modules
over an $A_\infty$-algebra $A$ are encoded either by 
comodule differentials on $BA \otimes M$ or twisted cochains $BA 
\to End(M)$. 

If $\tau$ is a generalized twisted cochain as above and $N$ is 
a DG comodule over $C$, we denote by  $N\otimes_\tau A$ the 
tensor product $N\otimes A$ with the differential 
$$
\delta = \delta_N\otimes 1 + 1 \otimes \delta_A
+\sum_{s \geq 2} (1 \otimes m_s)(1 \otimes \otimes 
\tau^{\otimes(s-1)} \otimes 1) (\Delta_N^{(s)} \otimes 1)
$$
where $m_s$ is the $s$-th product in $A$ and $\Delta_N^{(s)}:
N\to N\otimes \overline{C}^{\otimes (s-1)}$ is the iterated 
reduced coaction map. 
The infinite sum makes sense if $N$ is cocomplete, i.e. 
$N= \cup_i Ker (\Delta_N^{(i)})$. 

On the other hand, is $M$ is an $A_\infty$-module over $A$
with action maps
$m_s^M: M\otimes A^{\otimes (s-1)} \to M$ then denote by 
$M\otimes_\tau C$ the tensor product $M\otimes C$ with the differential
$$
\delta = \delta_M\otimes 1 + 1 \otimes \delta_C + 
\sum_{s \geq 2} (m_s^M\otimes 1)(1 \otimes \tau^{\otimes (s-1)}
\otimes 1)(1 \otimes \Delta^{(s)}).
$$
Observe that $N\otimes_\tau A$ is a quasi-free right 
$A_\infty$-module over $A$ and $M\otimes C$ is a quasi-free
DG comodule over $C$. See Section 3 in \cite{B2} on how 
to define the corresponding functors on morphisms, and other details.

\noindent
\textsl{Department of Mathematics, 103 MSTB}\\
\noindent
\textsl{University of California, Irvine}\\
\noindent
\textsl{Irvine, CA 92697, USA}\\
\noindent
\textsl{email: vbaranov@math.uci.edu}


\begin{thebibliography}{BKRS}

\bibitem[AAFR]{AAFR}

\'Alvarez, V.; Armario, J. A.; Frau, M. D.; Real, P.:
Transferring TTP-structures via contraction. 
\textit{Homology Homotopy Appl.} \textbf{7} (2005), no. 2, 41--54. 

\bibitem[B1]{B1} Baranovsky, V.: BGG correspondence for projective 
complete intersections,  \textit{Int. Math. Res. Not.}
\textbf{2005},  no. 45, 2759--2774.

\bibitem[B2]{B2} Baranovsky, V.: BGG correspondence for toric 
complete intersections, to appear in \textit{Moscow. Math. J.}. 

\bibitem[Br]{Br} 
Brown, R.: The twisted Eilenberg-Zilber theorem.  1965  
\textit{Simposio di Topologia (Messina, 1964)}.  pp. 33--37

%\bibitem[ABW]{ABW} 
%Akin, K.; Buchsbaum, D.; Weyman, J.: Schur functors and Schur complexes.  
%\textit{Adv. in Math.}  \textbf{44}  (1982), no. 3, 207--278. 



%\bibitem[CB]{CB} Crawley-Boevey, W.: 
%\textit{Lectures on representation theory and invariant theory}, notes
%available on author's homepage 
%http://www.amsta.leeds.ac.uk/\~\,pmtwc/

%\bibitem[F]{F} Fulton, W.: 
%Young tableaux. 
%With applications to representation theory and geometry. London Mathematical Society Student Texts, 35. \textit{Cambridge University Press, Cambridge}, 1997.

\bibitem[FHT]{FHT}F\'elix, Y.; Halperin, S.; Thomas, J.-C.:
Rational homotopy theory. Graduate Texts in Mathematics, 205. 
\textit{Springer-Verlag, New York}, 2001.



\bibitem[GLS]{GLS}
Gugenheim, V. K. A. M.; Lambe, L. A.; Stasheff, J. D.: Perturbation theory in differential homological algebra. II. 
\textit{Illinois J. Math.}  \textbf{35}  (1991),  no. 3, 357--373. 



\bibitem[K]{K}
Keller, B.: Introduction to $A$-infinity algebras and modules.  
\textit{Homology Homotopy Appl.}  \textbf{3}  (2001),  no. 1, 1--35.

\bibitem[KS]{KS} Kontsevich, M., Soibelman, Y.:
 Notes on A-infinity algebras, A-infinity categories and 
non-commutative geometry, I,  preprint math.RA/0606241.

\bibitem[LH]{LH}  Lef\`evre-Hasegawa, K.:  Sur les A-infini 
cat\'egories, 
preprint math.CT/0310337.

\bibitem[LM]{LM}
Lada, T.; Markl, M.: Strongly homotopy Lie algebras.  
\textit{Comm. Algebra}  \textbf{23}  (1995),  no. 6, 2147--2161.

\bibitem[LS]{LS}
Lambe, L.; Stasheff, J.: Applications of perturbation theory to 
iterated fibrations.  
\textit{Manuscripta Math.} \textbf{58}  (1987),  no. 3, 363--376.

\bibitem[M]{M} Merkulov, S. A.: 
Quantization of strongly homotopy Lie bialgebras, preprint 
math.QA/0612431. 

\bibitem[MS]{MS}
Markl, M.; Shnider, S.: Associahedra, cellular $W$-construction and 
products of $A_\infty$-algebras.  Trans. Amer. Math. Soc.  358  (2006),  
no. 6, 2353--2372.

\bibitem[MSS]{MSS}
Markl, M.; Shnider, S.; Stasheff, J.: Operads in algebra, topology and 
physics. Mathematical Surveys and Monographs, 96. \textit{American 
Mathematical 
Society, Providence, RI}, 2002.

\bibitem[PP]{PP}
Polishchuk, A.; Positselski, L.: Quadratic algebras. 
University Lecture Series, 37. 
\textit{American Mathematical Society, Providence, RI}, 2005.

\bibitem[S]{S} Sagan B.E.: The symmetric group. Representations, 
combinatorial algorithms, and symmetric functions. Second edition. 
Graduate Texts in Mathematics, 203. \textit{Springer-Verlag, New York}, 
2001.


\bibitem[SU]{SU}
Saneblidze, S.; Umble, R.: Diagonals on the permutahedra, multiplihedra 
and associahedra.  \textit{Homology Homotopy Appl.}  
\textbf{6}  (2004),  no. 1, 363--411.
\end{thebibliography}
\end{document}